\input{style/arxiv-ba.cfg}
\documentclass[ba,linksfromyear,preprint]{imsart}
\makeatletter
   \@ifpackageloaded{natbib}{}{\usepackage{natbib}}
\makeatother
\usepackage{latexsym,multirow,booktabs}
\usepackage{ifthen,pifont,comment,enumerate}

\pubyear{2015}
\volume{10}
\issue{2}
\firstpage{351}
\lastpage{377}
\doi{10.1214/14-BA916}

\makeatletter
\newcommand{\enquote}[1]{``#1''}
\def\ci{\perp\!\!\!\perp}
\def\diag{{\rm diag}}
\def\log{{\rm log}}

\def\tr{{\rm tr}}

\def\N{\small\textsc{N}}
\def\Exp{\small\textsc{Exp}}
\def\GIG{\small\textsc{GIG}}
\def\IW{\small\textsc{IW}}
\def\W{\small\textsc{W}}

\def\iG{\textit{G}}

\def\iGIW{\textit{GIW}}

\def\ci{\perp\!\!\!\perp}
\def\bA{\mathbf{A}}
\def\bB{\mathbf{B}}
\def\bC{\mathbf{C}}
\def\bD{\mathbf{D}}

\def\bI{\mathbf{I}}

\def\bS{\mathbf{S}}
\def\bs{\mathbf{s}}

\def\bu{\mathbf{u}}
\def\bV{\mathbf{V}}
\def\bv{\mathbf{v}}

\def\bw{\mathbf{w}}

\def\bY{\mathbf{Y}}
\def\by{\mathbf{y}}
\def\bZ{\mathbf{Z}}
\def\bz{\mathbf{z}}

\newcommand{\bfb}   {\mbox{\boldmath$\beta$}}
\newcommand{\bfg}   {\mbox{\boldmath$\gamma$}}

\newcommand{\bfomega}   {\mbox{\boldmath$\omega$}}
\newcommand{\bfOmega} {\mbox{\boldmath$\Omega$}}

\newcommand{\bfSigma} {\mbox{\boldmath$\Sigma$}}
\newcommand{\bfsigma} {\mbox{\boldmath$\sigma$}}

\newtheorem{proposition}{Proposition}

\makeatother

\begin{document}

\begin{frontmatter}
\title{Scaling It Up: Stochastic Search Structure Learning in Graphical Models}
\runtitle{Stochastic Search Structure Learning in Graphical Models}

\begin{aug}
\author[a]{\fnms{Hao} \snm{Wang}\corref{}\ead[label=e1]{haowang@sc.edu}}

\runauthor{Hao Wang}

\address[a]{Department of Epidemiology and Biostatistics, Michigan
State University, East Lansing, Michigan 48824, U.S.A., \printead{e1}}

\end{aug}

%
\begin{abstract}
Gaussian concentration graph models and covariance graph models are two
classes of graphical models that are useful for uncovering latent
dependence structures among multivariate variables. In the Bayesian
literature, graphs are often determined through the use of priors over
the space of positive definite matrices with fixed zeros, but these
methods present daunting computational burdens in large problems.
Motivated by the superior computational efficiency of continuous
shrinkage priors for regression analysis, we propose a new framework
for structure learning that is based on continuous spike and slab
priors and uses latent variables to identify graphs. We discuss model
specification, computation, and inference for both concentration and
covariance graph models. The new approach produces reliable estimates
of graphs and efficiently handles problems with
hundreds of variables.
\end{abstract}

%

\end{frontmatter}


\section{Introduction}
Graphical models use graph structures for modeling and making
statistical inferences regarding complex relationships among many variables.
Two types of commonly used graphs are undirected graphs, which
represent conditional dependence relationships among variables, and
bi-directed graphs, which encode marginal dependence among variables.
Structure learning refers to the problem of estimating unknown graphs
from the data and is usually carried out by sparsely estimating the
covariance matrix of the variables by assuming that the data follow a
multivariate Gaussian distribution. Under the Gaussian assumption,
undirected graphs are determined by zeros in the concentration matrix
and their structure learning problems are thus referred to as
concentration graph models; bi-directed graphs are determined by zeros
in the covariance matrix and their structure learning problems are thus
referred to as covariance graph models. This work concerns structure
learning in both concentration and covariance graph models.

Classical methods for inducing sparsity often rely on penalized
likelihood approaches \citep
{Banerjee2008,YuanLin2007,BienTib2011,Wang2012SC}. Model fitting then
uses deterministic optimization procedures such as coordinate descents.
Thresholding is another popular method for the sparse estimation of
covariance matrices for covariance graph models \citep
{BickelLevina2008,RothmanEtAl2009}; however, there is no guarantee that
the resulting estimator is always positive definite. Bayesian methods
for imposing sparsity require the specification of priors over the
space of positive definite matrices constrained by fixed zeros. Under
such priors, model determination is then carried out through stochastic
search algorithms to explore a discrete graphical model space. The
inherent probabilistic nature of the Bayesian framework permits
estimation via decision-theoretical
principles, addresses parameter and model uncertainty, and provides a
global characterization of the
parameter space. It also encourages the development of modular
structures that can be integrated
with more complex systems.

A major challenge in Bayesian graphical models is their computation.
Although substantial progress in computation for these two graphical
models has been made in recent years, scalability with dimensions
remains a significant issue, hindering the ability to adapt these
models to the growing demand for higher dimensional problems. Recently
published statistical papers on these two graphical models either focus
on small problems or report long computing times. In concentration
graph models, \citet{DLR2011} report that it takes approximately 1.5
days to fit a problem of 48 nodes on a dual-core 2.8 Ghz computer under
C; \citet{WangLi2012} report approximately two days for a 100 node
problem under MATLAB; and \citet{ChengLenkoski2012} report a computing
time of 1--20 seconds per one-edge update for a 150 node problem using
a 400 core server with 3.2 GHz CPU under R and C++. In covariance graph
models, \citet{SilvaGhahramani2009} fit problems up to $13$ nodes and
conclude that ``further improvements are necessary for larger problems.''

To scale up with dimension, this paper develops a new approach called
stochastic search structure learning (SSSL) to efficiently determine
covariance and concentration graph models. The central idea behind SSSL
is the use of continuous shrinkage priors characterized by binary
latent indicators in order to avoid the normalizing constant
approximation and to allow block updates of graphs. The use of
continuous shrinkage priors contrasts point-mass priors at zeros that
are used essentially by all existing methods for Bayesian structure
learning in these two graphical models.

The motivation for the SSSL comes from the successful developments of
continuous shrinkage priors in several related problems. In regression
analysis, continuous shrinkage priors were used in the seminal paper by
\citet{George93} in the form of a two component normal mixture for
selecting important predictors and these priors have recently garnered
substantial research attention as a computationally attractive
alternative for regularizing many regression coefficients (e.g.,
\citealt{ParkCasellla2008,GriffinBrown10,Armagan11}). In estimation of
covariance matrices, they are used for regularizing concentration
elements and have been shown to provide fast and accurate estimates of
covariance matrices \citep{Wang2011,KhondkerEtAl12}. In factor
analysis, they are used instead of point-mass priors \citep
{CarvalhoChang2008} for modeling factor loading matrices, efficiently
handling hundreds of variables \citep{BhattacharyaDunson2011}.

Nevertheless, the current work is fundamentally different from the
aforementioned works. The research focus here is the structure learning
of graphs, which is distinct from regression analysis, factor analysis,
and the pure covariance estimation problem that solely performs
parameter estimation without the structure learning of graphs. Although
continuous shrinkage priors generally perform very well in these
problems, little is known about their performance in problems of
structure learning. Because graphs are directly determined by
covariance matrices, the positive definiteness of any covariance
matrices poses methodological challenges to investigating prior
properties, as well as to the construction of efficient stochastic
search algorithms. The paper's main contributions are the development
and exploration of two classes of continuous shrinkage priors for
learning undirected and bi-directed graphs, as well as two efficient
block Gibbs samplers for carrying out the corresponding structure
learning tasks that fit problems of one or two hundred variables within
a few minutes.

\section{Background on graphical models}\label{sec:background}
Assume that $\by=(y_1,y_2,\ldots,y_p)^\prime$ is a
$p$-dimensional random vector following a multivariate normal
distribution $\N(0,\bfSigma)$ with mean of zero and covariance matrix
$\bfSigma\equiv(\sigma_{ij})$. Let $\bfOmega\equiv(\omega_{ij})
=\bfSigma^{-1}$ be the concentration matrix. Covariance and
concentration graph models are immediately related to $\bfSigma$ and
$\bfOmega$, respectively. Let $\bY$ be the $n\times p$ data matrix
consisting of $n$ independent samples of $\by$ and
let $\bS= \bY^\prime\bY$. The theory and existing methods for
structure learning are briefly reviewed in the next two sections.


\subsection{Concentration graph models}
Concentration graph models \citep{Dempster72} consider the
concentration matrix $\bfOmega$ and encode conditional dependence using
an undirected graph $G=(V,E)$, where $V =\{1,2,\ldots,p\}$ is a
non-empty set of vertices
and $E\subseteq\{(i,j): i<j \}$ is a set of edges representing
unordered pairs of vertices. The graph $G$ can also be indexed by a set
of $p(p-1)/2$ binary variables $\bZ=(z_{ij})_{i<j}$, where $z_{ij} = 1$
or $0$ according to whether edge $(i,j)$ belongs to $E$ and not.
Theoretically, the following properties are equivalent:
\[
z_{ij} = 0 \quad \Leftrightarrow\quad(i,j)\notin E \quad
\Leftrightarrow\quad y_i\ci y_j \mid\by_{-(ij)} \quad
\Leftrightarrow\quad \omega_{ij}=0,
\]
where $\by_{-(ij)}$ is the random vector containing all elements in $\by
$ except for $y_i$ and $y_j$, and ``$\Leftrightarrow$'' reads as ``if
and only if''.

In the Bayesian paradigm, concentration graph models are usually
modeled through hierarchical priors consisting of the following: \textit
{(i)} the conjugate \iG-Wishart prior $\bfOmega\sim\W_G(b,\bD)$ \citep
{dawid93,Roverato02} for $\bfOmega$ given the graph $\bZ$; and \textit
{(ii)} independent Bernoulli priors for each edge-inclusion indicator
$z_{ij}$ with inclusion probability $\pi$:
%
\begin{eqnarray}
p(\bfOmega\mid\bZ) &=& I_{GW}(b,\bD,\bZ)^{-1}
|\bfOmega|^{b-2\over2} \exp\{-{1\over2} \tr(\bD\bfOmega) \} 1_{\{
\bfOmega\in
M^+(\bZ)\} }, \label{eq:GW} \\
p( \bZ) &=& \prod_{i<j}\big\{ \pi^{z_{ij}} (1-\pi)^{1-z_{ij}} \big\},
\label{eq:ind:bern}
\end{eqnarray}
where $b$ is the degrees-of-freedom parameter, $\bD$ is the location
parameter, $I_{GW}(b,\bD,\bZ)$ is the normalizing constant, and $M^+(\bZ
)$ is the cone of symmetric positive definite matrices with
off-diagonal entries $\omega_{ij}=0$
whenever $z_{ij}=0$. As for the hyperparameters, common choices are
$b=3,\bD=\bI_p$ with $\bI_p$ being the $p\times p$ identity matrix and
$\pi=2/(p-1)$ \citep{Jones05}. Under (\ref{eq:GW})--(\ref
{eq:ind:bern}), some methods (e.g.,
\citealt{Jones05,ScottCarvalho08,LenkoskiDobraJCGS}) learn $\bZ$
directly through its posterior distribution over the model space $p(\bZ
\mid\bY) \propto p(\bY\mid\bZ)p(\bZ) $. Other methods learn $\bZ$
indirectly through sampling over the joint space of graphs and
concentration matrices $p(\bfOmega,\bZ\mid\bY)$ \citep
{GiudiciGreen99Biom,DLR2011,WangLi2012}. Regardless of the types of
algorithms, two shared features of these methods cause the framework
(\ref{eq:GW})--(\ref{eq:ind:bern}) to be inefficient for larger $p$
problems. The first of these features is that graphs are updated in a
one-edge-at-a-time manner, meaning that sweeping through all possible
edges requires a loop of $O(p^2)$ iterations. The second feature is
that the normalizing constant $I_{GW}(b,\bD,\bZ)$ for
non-decomposable graphs requires approximation. The commonly used Monte
Carlo approximation proposed by \citet{Massam05} is not only unstable
in some situations but also requires a matrix completion step of time complexity
$O(M p^4)$ for $M$ Monte Carlo samples, making these methods
unacceptably slow
in large graphs. Recent works by \citet{WangLi2012} and \citet
{ChengLenkoski2012} propose the use of exchange algorithms to avoid the
Monte Carlo approximation. However, the computational burden remains
daunting; empirical experiments in these papers suggest it would take
several days to complete the fitting for problems of $p\approx100$ on
a desktop computer.

In the classical formulation, concentration graphs are induced by
imposing a graphical lasso penalty on $\bfOmega$ in order to encourage
zeros in the penalized maximum likelihood estimates of $\bfOmega$
(e.g., \citealt{YuanLin2007,Friedman08,Rothman2008}). In particular,
the standard graphical lasso problem
is to maximize the penalized log-likelihood
\[
\log(\det\bfOmega)-\tr({\bS\over n}\bfOmega) - \rho
||\bfOmega||_1,
\]
over the space of positive definite matrices $M^+$, with $\rho\geq0$
as the shrinkage parameter and $||\bfOmega||_1 =
\sum_{1\leq i,j \leq p} |\omega_{ij}|$ as the $L_1$-norm of
$\bfOmega$. The graphical lasso problem has a Bayesian interpretation
\citep{Wang2011}. Its estimator is equivalent to the maximum a
posteriori estimation
under the following prior for $\bfOmega$:
%
\begin{eqnarray} \label{eq:bglasso}
p(\bfOmega) &=& C^{-1} \prod_{1 \leq i,j \leq p} \big\{ \exp(-\lambda
|\omega_{ij}| )\big\} 1_{\bfOmega\in M^+}, \label{eq:bglasso}
\end{eqnarray}
where $C$ is the normalizing constant. By exploiting the scale mixture
of normal representation, \citet{Wang2011} shows that fitting (\ref
{eq:bglasso}) is very efficient using block Gibbs samplers for up to
the lower hundreds of variables.

A comparison between the two Bayesian methods (\ref{eq:GW})--(\ref
{eq:ind:bern}) and (\ref{eq:bglasso}) helps to explain the intuition
behind the proposed SSSL. Model (\ref{eq:GW})--(\ref{eq:ind:bern})
explicitly treats a graph $\bZ$ as an unknown parameter and considers
its posterior distribution, which leads to straightforward Bayesian
inferences about graphs. However, it is slow to run due to the
one-edge-at-a-time updating and the normalizing constant approximation.
In contrast, Model (\ref{eq:bglasso}) uses continuous priors, enabling
a fast block Gibbs sampler that updates $\bfOmega$ one column at a time
and avoids normalizing constant evaluation. However, no graphs are used
in the formulation, and thus this approach does not constitute a formal
Bayesian treatment of structure learning. Still, a better approach
might be developed by using the best aspects of the two methods. That
is, such a method would allow explicit structure learning, as in (\ref
{eq:GW})--(\ref{eq:ind:bern}), while maintaining good scalability, as
in (\ref{eq:bglasso}). This possibility is exactly the key of SSSL.
Similar insights also apply to the covariance graph models described below.

\subsection{Covariance graph models}
Covariance graph models \citep{CoxWermuth1993} consider the covariance
matrix $\bfSigma$ and encode the marginal dependence using a
bi-directed graph $G=(V,E)$, where each edge has bi-directed
arrows instead of the full line used by an undirected graph. Similar to
concentration graph models, the covariance graph $G$ can also be
indexed by binary variables $\bZ=(z_{ij})_{i<j}$. Theoretically, the
following properties are equivalent:
\[
z_{ij} = 0 \quad\Leftrightarrow \quad(i,j)\notin E \quad
\Leftrightarrow \quad y_i\ci y_j \quad \Leftrightarrow \quad \sigma_{ij}=0.
\]
In the Bayesian framework, structure learning again relies on the
general hierarchical priors $p(\bfSigma,\bZ)=p(\bfOmega\mid\bZ)p(\bZ
)$. For $p(\bfSigma\mid\bZ)$, \citet{SilvaGhahramani2009} propose the
conjugate \iG-inverse Wishart prior $\bfSigma\sim\IW_{G}(b,\bD)$ with
the density as:
%
\begin{eqnarray}\label{eq:GIW} p(\bfSigma\mid\bZ) &=& I^{-1}_{\iGIW
}(b,\bD,\bZ)|\bfSigma|^{-{b+2p\over2}}\exp\{-{1\over2}\tr(\bD\bfSigma
^{-1})\}1_{\bfSigma\in M^+(\bZ)},
\end{eqnarray}
where $b$ specifies the degrees of freedom, $\bD$ is the location
parameter, and $I_{GIW}(b,\bD,\bZ)$ is the normalizing constant.
Structure learning is then carried out through the marginal likelihood
function $p(\bY\mid\bZ) = (2\pi)^{-{np/ 2}}{ I_\iGIW(b+n,\bD+\bS,\bZ
)/ I_\iGIW(b,\bD,\bZ)}$.
Unfortunately, the key quantity of the normalizing constant $I_\iGIW
(b,\bD,\bZ)$ is analytically intractable, even for decomposable graphs.
\citet{SilvaGhahramani2009} propose a Monte Carlo approximation via an
importance sampling algorithm, which becomes computationally infeasible
for $p$ beyond a few dozen. Their empirical experiments are thus
limited to small problems (i.e., $p<20$). Later, \citet{Bala11}
investigate a broad class of priors for decomposable covariance graph
models that embed (\ref{eq:GIW}) as a special case. They also derive
closed-form normalizing constants for decomposable homogeneous graphs
which account for only a tiny portion of the overall graph space.
Despite these advances, the important question of scalability to higher
dimensional problems remains almost untouched.

In the classical framework, the earlier literature focuses on maximum
likelihood estimates and likelihood ratio test procedures (e.g.,
\citealt{Kauermann1996,Wermuth2006,Chaudhuri2007}). Later, two general
types of approaches are proposed to estimate zeros in the covariance
elements. The first is the thresholding procedure, which sets $\sigma
_{ij}$ to be zero if its sample estimate is below a threshold \citep
{BickelLevina2008,RothmanEtAl2009,CaiLiu2011}. Another approach is
motivated by the lasso-type procedures. \citet{BienTib2011} propose a
covariance graphical lasso procedure for simultaneously estimating
covariance matrix and marginal dependence structures. Their method is
to minimize the following objective function:
%
\begin{eqnarray} \label{eq:obj}
\log(\det\bfSigma)+\tr({\bS\over n}\bfSigma^{-1}) + \rho||\bfSigma||_1,
\end{eqnarray}
over the space of positive definite matrices $M^+$, with $\rho\geq0$
as the shrinkage parameter. In comparison with thresholding and
likelihood ratio testing methods, this approach has the advantage of
guaranteeing the positive definiteness of the estimated $\bfSigma$.
Although a Bayesian version of (\ref{eq:obj}) has not been explored
previously, its derivation is straightforward through the prior
%
\begin{eqnarray}\label{eq:bcglasso}
p(\bfSigma) &=& C^{-1} \prod_{1 \leq i,j \leq p} \big\{ \exp(-\lambda
|\sigma_{ij}| )\big\} 1_{\bfSigma\in M^+},
\end{eqnarray}
In light of the excellent performance of the Bayesian concentration
graphical lasso (\ref{eq:bglasso}) reported in \citet{Wang2011}, we
hypothesize that (\ref{eq:bcglasso}) shares similar performances. In
fact, we have developed a block Gibbs sampler for (\ref{eq:bcglasso})
and found that it gives a shrinkage estimation of $\bfSigma$ and is
computationally efficient, although it provides no explicit treatment
of the graph $\bZ$. The detailed algorithm and results are not reported
in this paper but are available upon request. Comparing (\ref{eq:GIW})
and (\ref{eq:bcglasso}) suggests that again, the different strengths of
(\ref{eq:GIW}) and (\ref{eq:bcglasso}) might be combined to provide a
better approach for structure learning in covariance graph models.

\section{Continuous spike and slab priors for positive definite
matrices}\label{sec:cov}
Let $\bA=(a_{ij})_{p\times p}$ denote a $p$-dimensional covariance or
concentration matrix; that is, $\bA=\bfSigma$ or $\bfOmega$. SSSL uses
the following new prior for $\bA$:
%
\begin{equation}\label{eq:prior:pd}
p(\bA) = \{C(\theta)\}^{-1} \prod_{i<j} \bigg\{ (1-\pi) \N(a_{ij} \mid
0,v_0^2) + \pi\N(a_{ij}\mid0,v_1^2) \bigg\} \prod_i \Exp(a_{ii}\mid
{\lambda\over2}) 1_{\{\bA\in M^+\}},
\end{equation}
where $\textsc{N}(a\mid0,v^2) $ is the density
function of a normal random variable
with mean $0$ and variance $v^2$ evaluated at $a$, $\textsc{Exp}(a\mid
\lambda)$ represents the exponential density function of the form $p(a)
= \lambda\exp(-\lambda x)1_{a>0}$, and $1_{\{\cdot\}}$ is the
indicator function. The parameter spaces are $v_0>0$, $v_1>0$, $\lambda
>0$ and $\pi\in(0,1)$, and the set of all parameters is denoted as
$\theta= \{v_0,v_1,\pi,\lambda\}$. The values of $v_0$ and $v_1$ are
further set to be small and large, respectively. The term $C(\theta)$
represents the normalizing constant that ensures the integration of the
density function $p(\bA)$ over the space $M^+$ is one, and it depends
on $\theta$. The first product symbol multiplies $p(p-1)/2$ terms of
two-component normal mixture densities for the off-diagonal elements,
connecting this prior to the classical and Bayesian graphical lasso
methods through the familiar framework of normal mixture priors for
$a_{ij}$. The second product symbol multiplies $p$ terms of exponential
densities for the diagonal elements. The two-component normal mixture
density plays a critical role in structure learning, as will be clear below.

Prior (\ref{eq:prior:pd}) can be defined by introducing binary latent
variables $\bZ\equiv(z_{ij})_{i<j}\in\mathcal{Z} \equiv\{0,1\}
^{p(p-1)/2}$ and a hierarchical
model:\vadjust{\eject}
\begin{align}
p(\bA\mid\bZ,\theta) & = \{C({\bZ,v_0,v_1,\lambda})\}^{-1} \prod
_{i<j} \N(a_{ij}\mid0,v^2_{z_{ij}}) \prod_i \Exp(a_{ii}\mid{\lambda
\over2}), \label{eq:hier:1} \\
p( \bZ\mid\theta) & = \{C({\theta})\}^{-1} C({\bZ,v_0,v_1,\lambda})
\prod_{i<j}\big\{ \pi^{z_{ij}} (1-\pi)^{1-z_{ij}} \big\}, \label{eq:hier:2}
\end{align}
where $v_{z_{ij}} = v_0$ or $v_1$ if $z_{ij}=0$ or $1$.
The intricacy here is the two terms of $C({\bZ,v_0,v_1,\lambda})$. Note
that $C({\bZ,v_0,v_1,\lambda})\in(0,1)$ because it is equal to the
integration of the product of normal densities and exponential
densities over a constrained space $M^+$. Thus, (\ref{eq:hier:1}) and
(\ref{eq:hier:2}) are proper distributions. The joint distribution of
$(\bA,\bZ)$ acts to cancel out the two terms of $C({\bZ,v_0,v_1,\lambda
})$ and results in a marginal distribution of $\bA$, as in (\ref{eq:prior:pd}).


The rationale behind using $\bZ$ for structure learning is as follows.
For an appropriately chosen small value of $v_0$, the event $z_{ij}=0$
means that $a_{ij}$ comes from the concentrated component $\N
(0,v_0^2)$, and so $a_{ij}$ is likely to be close to zero and can
reasonably be estimated as zero. For an appropriately chosen large
value of $v_1$, the event $z_{ij}=1$ means that $a_{ij}$ comes from the
diffuse component $\N(0,v_1^2)$ and so $a_{ij}$ can be estimated to be
substantially different from zero. Because zeros in $\bA$ determine
missing edges in graphs, the latent binary variables $\bZ$ can be
viewed as edge-inclusion indicators. Given data $\bY$, the posterior
distribution of $\bZ$ provides information about graphical model
structures. The remaining questions are then how to specify parameters
$\theta$ and how to perform posterior computations.

\subsection{Choice of $\pi$}
From (\ref{eq:hier:2}), the hyperparameter $\pi$ controls the prior
distribution of the edge-inclusion indicators in $\bZ$. The choice of
$\pi$ should thus reflect the prior belief about what the graphs will
be in the final model. In practice, such prior information is often
summarized via the marginal prior edge-inclusion probability $\Pr
(z_{ij}=1)$. Specifically, a prior for $\bZ$ is chosen such that the
implied edge-inclusion probability of edge $(i,j)$ meets the prior
belief about the chance of the existence of edge $(i,j)$. For example,
the common choice $\Pr(z_{ij}=1) = 2/(p-1)$ reflects the prior belief
that the expected number of edges is approximately ${p \choose2}\Pr
(z_{ij}=1) = p$. Another important reason that $\Pr(z_{ij}=1)$ is used
for calibrating priors over $\bZ$ is that the posterior inference about
$\bZ$ is usually based upon the marginal posterior probability of $\Pr
(z_{ij}=1 \mid\bY)$. For example, the median probability graph, the
graph consisting of those edges whose marginal edge-inclusion
probability exceeds $0.5$, is often used to estimate $G$ \citep
{Wang2010}. Focusing on the marginal edge-inclusion probability allows
us to understand how the posterior truly responds to the data.

Calibrating $\pi$ according to $\Pr(z_{ij}=1)$ requires knowledge of
the relation between $\Pr(z_{ij}=1)$ and $\pi$. From (\ref{eq:hier:2}),
the explicit form of $\Pr(z_{ij}=1)$ as a function of $\pi$ is
unavailable because of the intractable term $C({\bZ,v_0,v_1,\lambda})$.
A comparison between (\ref{eq:hier:2}) and (\ref{eq:ind:bern}) helps
illustrate the issue. Removing $C({\bZ,v_0,v_1,\lambda})$ from (\ref
{eq:hier:2}) turns it into (\ref{eq:ind:bern}) but will not cancel out
the term $C({\bZ,v_0,v_1,\lambda})$ in (\ref{eq:hier:1}) for the
posterior distribution of $\bZ$. Tedious and unstable numerical
integration is then necessary to evaluate $C({\bZ,v_0,v_1,\lambda})$ at
each iteration of sampling $\bZ$. Inserting $C({\bZ,v_0,v_1,\lambda})$
into (\ref{eq:hier:2}) cancels $C({\bZ,v_0,v_1,\lambda})$ in (\ref
{eq:hier:1}) in\vadjust{\goodbreak} the posterior, thus facilitating computation, yet
concerns might be raised about such a ``fortunate'' cancelation. For
example, \citet{Murray07} notes that a prior that cancels out an
intractable normalizing constant in the likelihood would depend on the
number of data points and would also be so extreme that it would dominate
posterior inferences. These two concerns appear to be not problematic
in our case. The prior (\ref{eq:hier:2}) is for the hyperparameter $\bZ
$, rather than for the parameter directly involved in the likelihood;
thus it does not depend on sample size. Instead, the prior also only
introduces mild bias without dominating the inferences, as shown below.

%
%
\begin{figure}[b!]
\centering
\includegraphics{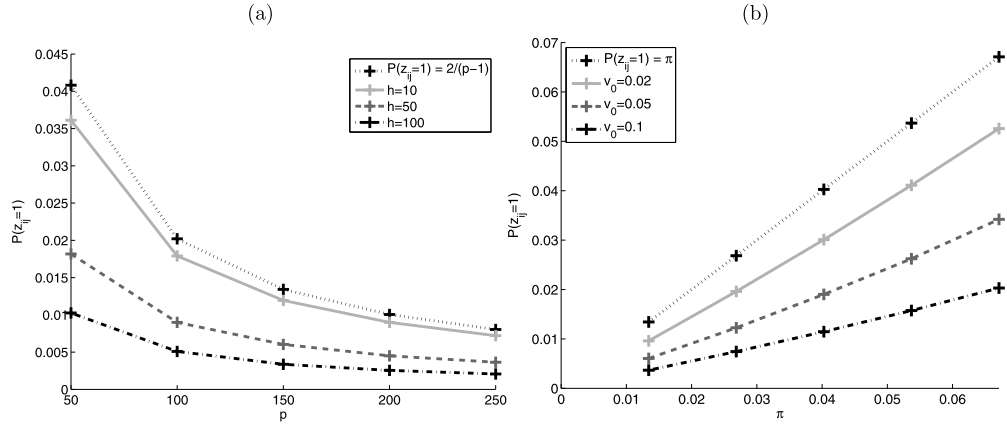}
\caption{The implied prior marginal edge-inclusion probability $\Pr
(z_{ij}=1)$ from the prior (\ref{eq:hier:1})--(\ref{eq:hier:2}) as a
function of $p$ at different $h$ (left) and a function of $\pi$ at
different $v_0$ (right), together with two reference curves of $\Pr
(z_{ij}=1)=2/(p-1)$ (left) and of $\Pr(z_{ij}=1)=\pi$ (right). }
\label{fig:prior:z}
\end{figure}

To investigate whether the cancelation of $C({\bZ,v_0,v_1,\lambda})$
causes the prior to be too extreme, we compute $\Pr(z_{ij}=1)$
numerically from Monte Carlo samples generated by the algorithm in
Section \ref{sec:con}. In (\ref{eq:hier:1})--(\ref{eq:hier:2}), we
first fix $\pi= 2/(p-1), v_0 = 0.05$, and $\lambda=1$, and then vary
the dimension $p \in\{50,100,150,200,250\}$ and $v_1=hv_0$ with $h \in
\{10, 50, 100\}$. Panel (a) of Figure \ref{fig:prior:z} displays these
estimated $\Pr(z_{ij}=1)$ as a function of $p$ for different $h$
values. As a reference, the curve $\Pr(z_{ij}=1) = 2/(p-1)$ is also
plotted. The most noticeable pattern is that all three curves
representing the implied $\Pr(z_{ij}=1)$ from (\ref{eq:hier:2}) are
below the reference curve, suggesting that there is a downward bias
introduced by $C({\bZ,v_0,v_1,\lambda})$ on $\Pr(z_{ij}=1)$. The bias
is introduced by the fact that the positive definite constraint on $\bA
$ favors a small $v_0$, specified by $z_{ij}=0$, over a large $v_1$,
specified by $z_{ij}=1$. We also see that the bias is larger at larger
values of $h$, at which the impact of positive definite constraints is
more significant. Next, we fix $p=150, h=50,$ and $\lambda=1$ and vary
$v_0 \in\{0.02,0.05,0.1\}$ and $\pi\in\{
2/(p-1),4/(p-1),6/(p-1),8/(p-1),10/(p-1)\}$. Panel (b) of Figure \ref
{fig:prior:z} displays these implied $\Pr(z_{ij}=1)$ as a function of
$\pi$ for different $v_0$ values. Again, as a reference, the curve $\Pr
(z_{ij}=1) = \pi$ is plotted. The downward bias of $\Pr(z_{ij}=1)$
relative to $\pi$ continues to exist and is larger at larger values of
$v_0$ or $\pi$ because the positive definite constraint on $\bA$ forces
$z_{ij}=0$ to be chosen more often when $v_0$ or $\pi$ is large.
Nevertheless, the downward bias seems to be gentle, as $\Pr(z_{ij}=1)$
is never extremely small; consequently the prior (\ref
{eq:hier:1})--(\ref{eq:hier:2}) is able to let the data reflect the $\bZ
$ if the likelihood is strong.




Another concern is that the lack of analytical relation between $\Pr
(z_{ij}=1)$ and $\pi$ might raise challenges against the incorporation
of prior information about $\Pr(z_{ij}=1)$ into $\pi$. This problem can
be side-stepped by prior simulation and interpolation. Take a $p=150$
node problem as an example. If the popular choice $\Pr
(z_{ij}=1)=2/(p-1)=0.013$ is desirable, interpolating the points $(\pi
,\Pr(z_{ij}=1))$ in Panel (b) of Figure \ref{fig:prior:z} suggests that
$\pi$ should be set approximately at $0.018$, $0.027$, and $0.048$ for
$v_0 =0.02$, $0.05$, and $0.1$, respectively. Our view is that
obtaining these points under prior (\ref{eq:hier:1})--(\ref{eq:hier:2})
is much faster than evaluating $C({\bZ,v_0,v_1,\lambda})$ at each
configuration of $\bZ$ under the traditional prior (\ref{eq:ind:bern}).

%
\begin{figure}[t!]
\centering
\includegraphics{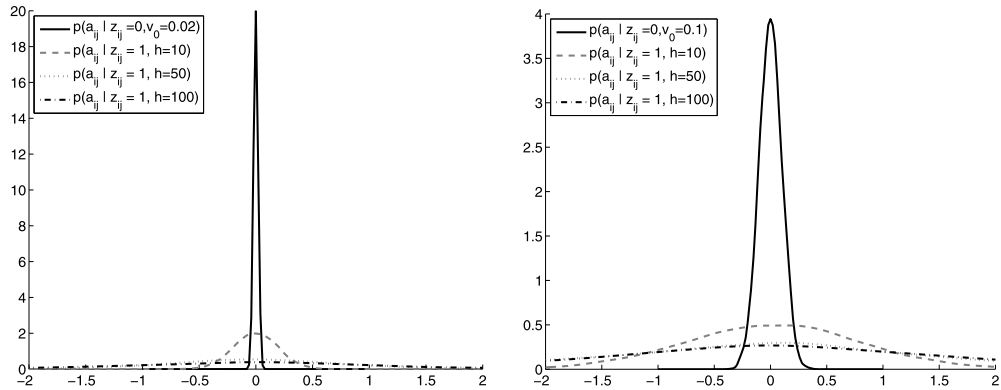}
\caption{The univariate density of $p(a_{ij} \mid z_{ij})$ for
different values of $h$ and $v_0$.}
\label{fig:prior:a}\vspace*{-3pt}
\end{figure}

\subsection{Choice of $v_0$ and $v_1$}
From (\ref{eq:hier:1}), the choice of $v_0$ should be such that if the
data supports $z_{ij}=0$ over $z_{ij}=1$, then $a_{ij}$ is small enough
to be replaced by zero. The choice of $v_1$ should be such that, if the
data favor $z_{ij}=1$ against $z_{ij}=0$, then $a_{ij}$ can be
accurately estimated to be substantially different from zero. One
general strategy for choosing $v_0$ and $v_1$, as recommended by \citet
{George93}, is based on the concept of practical significance.
Specifically, suppose a small $\delta$ can be chosen such that
$|a_{ij}|<\delta$ might be regarded as practically insignificantly
different from zero. Incorporating such a prior belief is then achieved
by choosing $v_0$ and $v_1$, such that the density $p(a_{ij} \mid
z_{ij}=0)$ is larger than the density $p(a_{ij} \mid z_{ij}=1)$,
precisely within the interval $(-\delta,\delta)$. An explicit
expression of $v_0$ as a function of $\delta$ and $h$ can be derived
when $p(a_{ij} \mid z_{ij})$ are normal. However, the implied
distribution $p(a_{ij} \mid z_{ij})$ from (\ref{eq:hier:1})--(\ref
{eq:hier:2}) is neither normal nor even analytically tractable. A
numerical study will illustrate some aspects of $p(a_{ij} \mid z_{ij})$.

Figure \ref{fig:prior:a} draws the Monte Carlo estimated density of
$p(a_{ij} \mid z_{ij}=0)$ and $p(a_{ij} \mid z_{ij}=1)$ for
several\vadjust{\goodbreak}
settings of $v_0$ and $h$. In all cases, there is a clear separation
between $p(a_{ij} \mid z_{ij}=0)$ and $p(a_{ij} \mid z_{ij}=1)$, with a
larger $h$ resulting in a sharper separation. This property of
separation between a small and a large variance component is the
essence of the prior for structural learning that aims to separate
small and large $a_{ij}$s. Clearly, the marginal densities are no
longer normal. For example, the density of $p(a_{ij} \mid z_{ij}=0)$ is
more spiky than that of $\N(0,v_0)$; the difference between $p(a_{ij}
\mid z_{ij}=1)$ when $h=50$ and $h=100$ is less clear than the
difference between $\N(0,2500v_0^2)$ and $\N(0,10000v_0^2)$. The lack
of an explicit form of $p(a_{ij} \mid z_{ij})$ makes the strategies of
analytically calculating $v_0$ from the threshold $\delta$ infeasible.
Numerical methods that estimate $p(a_{ij} \mid z_{ij})$ from Markov
chain Monte Carlo (MCMC) samples might be used to choose $v_0$
according to $\delta$ from a range of possible values.

Another perspective is that, when $v_0$ is chosen to be very small and
$h$ is chosen to be very large, the prior for $a_{ij}$ is close to a
point-mass mixture that selects any $a_{ij}\neq0$ as an edge. Because
the point-mass mixture prior provides a noninformative method of
structure learning when the threshold $\delta$ cannot be meaningfully
specified, it makes sense to choose $v_0$ to be as small as possible,
but not so small that it could cause MCMC convergence issues, and to
choose $v_1$ to allow for reasonable values of $a_{ij}$. In our
experiments with standardized data, MCMC converges quickly and mixes
quite well, as long as $v_0 \geq0.01$ and $h\leq1000$.

\subsection{Choice of $\lambda$}
The value of $\lambda$ controls the distribution of the diagonal
elements $a_{ii}$. Because the data are usually standardized, a choice
of $\lambda=1$ assigns substantial probability to the entire region of plausible
values of $a_{ii}$, without overconcentration or overdispersion. From
our experience, the data generally contain sufficient information about
the diagonal elements, and the structure learning results are
insensitive to a range of $\lambda$ values, such as $\lambda=5$ and $10$.

\section{Fast block Gibbs samplers}\label{sec:comp}
The primary advantage of the SSSL prior (\ref{eq:hier:1})--(\ref
{eq:hier:2}) over traditional approaches is its scalability to larger
$p$ problems. The reduction in computing time comes from two
improvements. One is that (\ref{eq:hier:1})--(\ref{eq:hier:2}) enable
block updates of all $p(p-1)/2$ edge-inclusion indicators in $\bZ$
simultaneously, while other methods only update one edge-inclusion
indicator $z_{ij}$ at a time. The other is that there is no need of a
Monte Carlo approximation of the intractable constants, while all of
the other methods require some sort of Monte Carlo integration to
evaluate any graph. The general sampling scheme for generating
posterior samples of graphs is to sample from the joint distribution
$p(\bA,\bZ\mid\bY)$ by iteratively generating from $p(\bA\mid\bZ,\bY
)$ and $p(\bZ\mid\bA,\bY)$. The first conditional $p(\bA\mid\bZ,\bY
)$ is sampled in a column-wise manner and the second conditional $p(\bZ
\mid\bA,\bY)$ is generated all at once. The details depend on whether
$\bA= \bfOmega$ for concentration graph models or $\bA= \bfSigma$ for
covariance graph models, and they are described below.

\subsection{Block Gibbs samplers for concentration graph models}\label{sec:con}
Consider concentration graph models with $\bA= \bfOmega$ in the
hierarchical prior (\ref{eq:hier:1})--(\ref{eq:hier:2}).
To sample from $p(\bfOmega\mid\bZ,\bY)$, the following proposition
provides necessary conditional distributions. The proof is in the Appendix.
\begin{proposition}\label{propsn1}
Suppose $\bA= \bfOmega$ in the hierarchical prior (\ref
{eq:hier:1})--(\ref{eq:hier:2}). Focus on the last column and row. Let
$\bV= (v_{z_{ij}}^2)$ be a $p\times p$ symmetric matrix with zeros in
the diagonal entries and $(v_{ij}^2)_{i<j}$ in the upper diagonal
entries. Partition $\bfOmega, \bS=\bY'\bY$ and $\bV$ as follows:
%
\begin{equation}\label{eq:partition:con}
\bfOmega= \left(
\begin{array}{cc}
\bfOmega_{11},\bfomega_{12} \\ \bfomega_{12}^\prime,\omega_{22}
\end{array}
\right),\quad\bS= \left(
\begin{array}{cc}
\bS_{11},\bs_{12} \\ \bs_{12}^\prime,s_{22}
\end{array}
\right), \quad\bV= \left(
\begin{array}{cc}
\bV_{11},\bv_{12} \\ \bv_{12}^\prime,0
\end{array}
\right).
\end{equation}
%
Consider a change of variables: $(\bfomega_{12},\omega_{22})
\rightarrow(\bu=\bfomega_{12}, v = \omega_{22} - \bfomega_{12}^\prime
\bfOmega_{11}^{-1} \bfomega_{12})$. We then have the full conditionals:
\begin{eqnarray}\label{eq:con:cond}
(\bu\mid-) \sim\N(-\bC\bs_{12},\bC), \quad\textrm{and } (v \mid-)
\sim\textsc{Ga}({n\over2}+1,{s_{22} + \lambda\over2}),
\end{eqnarray}
where $\bC= \{(s_{22}+\lambda)\bfOmega_{11}^{-1}+ \diag(\bv
_{12}^{-1})\}^{-1}$.
\end{proposition}
Permuting any column to be updated to the last one and using \eqref
{eq:con:cond} will lead to a simple block Gibbs step for generating
$(\bfOmega\mid\bZ,\bY)$. For $p(\bZ\mid\bfOmega,\bY)$, prior (\ref
{eq:hier:1})--(\ref{eq:hier:2}) implies all $z_{ij}$ are independent
Bernoulli with probability
%
\begin{eqnarray}\label{eq:con:z}
\Pr(z_{ij} &=& 1 \mid\bfOmega,\bY) = {\N(\omega_{ij}\mid0,v^2_1) \pi
\over\N(\omega_{ij}\mid0,v^2_1) \pi+ \N(\omega_{ij}\mid0,v_0^2)
(1-\pi) }.
\end{eqnarray}

A closer look at the conditional distributions of the last column $\bu
=\bfomega_{12}$ in \eqref{eq:con:cond} and the corresponding
edge-inclusion indicator vector $\bfg\equiv(\gamma_1,\ldots,\gamma
_{p-1})'= (z_{1p},\ldots,z_{p-1,p})'$ in \eqref{eq:con:z} reveals
something interesting. These distributions look like the Gibbs samplers
used in the stochastic search variable selection (SSVS) algorithm \citep
{George93}. Indeed, consider $\bfb\equiv(\beta_1,\ldots,\beta
_{p-1})'= -\bu$ and note that $s_{22}=n$ for standardized data. If
$\bfOmega_{11}^{-1}= {1\over n} \bS_{11}$ and $\lambda=0$,
then (\ref{eq:con:cond}) implies that
\[
(\bfb\mid\bz_{12},\bY) \sim\N\bigg[\big\{\bS_{11}+\diag(\bv
_{12}^{-1})\big\}^{-1} \bs_{12}, \big\{\bS_{11}+ \diag(\bv
_{12}^{-1})\big\}^{-1} \bigg],
\]
and \eqref{eq:con:z} implies that
\[
\Pr(\gamma_{j} = 1 \mid\bfb) = {\N(\beta_{j}\mid0,v_1^2) \pi\over\N
(\beta_{j}\mid0, v_1^2) \pi+ \N(\beta_{j}\mid0,v_0^2) (1-\pi) },
\quad j=1,\ldots,p-1,
\]
which are exactly the Gibbs sampler of SSVS for the $p$-th variable.
Thus, the problem of SSSL for concentration graph models can be viewed
as a $p$-coupled SSVS regression problem, as the use of $\bfOmega
_{11}^{-1}$ in the conditional distribution of $\bfomega_{12}$ in place
of $\bS_{11}$ shares information across $p$ regressions in a coherent
fashion. This interesting connection has not been noted elsewhere, to
the best of our knowledge.

\subsection{Block Gibbs samplers for covariance graph models}\label{sec:cov}
Now, consider covariance graph models with $\bA=\bfSigma$ in the
hierarchical prior (\ref{eq:hier:1})--(\ref{eq:hier:2}). To sample from
$p(\bfSigma\mid\bZ,\bY)$, the following proposition provides
necessary conditional distributions; its proof is in the Appendix.
\begin{proposition}\label{propsn2}
Suppose $\bA= \bfSigma$ in the hierarchical prior (\ref
{eq:hier:1})--(\ref{eq:hier:2}). Focus on the last column and row. Let
$\bV= (v_{z_{ij}}^2)$ be a $p\times p$ symmetric matrix with zeros in
the diagonal entries and $(v_{z_{ij}}^2)_{i<j}$ in the upper diagonal
entries. Partition $\bfSigma, \bS$ and $\bV$ as follows:
%
\begin{equation}\label{eq:partition:cov}
\bfSigma= \left(
\begin{array}{cc}
\bfSigma_{11},\bfsigma_{12} \\ \bfsigma_{12}^\prime,\sigma_{22}
\end{array}
\right),\quad\bS= \left(
\begin{array}{cc}
\bS_{11},\bs_{12} \\ \bs_{12}^\prime,s_{22}
\end{array}
\right), \quad\bV= \left(
\begin{array}{cc}
\bV_{11},\bv_{12} \\ \bv_{12}^\prime,0
\end{array}
\right).
\end{equation}
Consider a change of variables: $(\bfsigma_{12},\sigma_{22})\rightarrow
(\bu=\bfsigma_{12}, v = \sigma_{22} - \bfsigma_{12}^\prime\bfSigma
_{11}^{-1} \bfsigma_{12}).$ We then have the full conditionals:
\begin{align}\label{eq:beta:gamma} (\bu\mid\bY,\bZ,\bfSigma_{11},v) &
\sim \N\bigg[ \big\{\bB+ \diag(\bv_{12}^{-1})\big\}^{-1} \bw, \{\bB
+\diag(\bv_{12}^{-1})\big\}^{-1} \bigg], \nonumber\\ (v \mid\bY,\bZ
,\bfSigma_{11},\bu) & \sim\GIG(1-n/2,\lambda, \bu^\prime\bfSigma
_{11}^{-1} \bS_{11} \bfSigma_{11}^{-1} \bu-2\bs_{12}^\prime\bfSigma
_{11}^{-1}\bu+ s_{22}),
\end{align}
where $\bB= \bfSigma_{11}^{-1} \bS_{11} \bfSigma_{11}^{-1}v^{-1}+\lambda
\bfSigma_{11}^{-1}$, $\bw= \bfSigma_{11}^{-1}\bs_{12}v^{-1}$, and $\GIG
(q,a,b)$ denotes the generalized inverse Gaussian distribution with a
probability density function:
\[
p(x) = \frac{(a/b)^{q/2}}{2 K_q(\sqrt{ab})} x^{(p-1)} e^{-(ax + b/x)/2},
\]
with $K_q$ as a modified Bessel function of the second kind.
\end{proposition}
Surprisingly, Proposition (\ref{propsn2}) shows that the conditional
distribution of any column (row) in $\bfSigma$ is also multivariate
normal. This suggests direct column-wise block Gibbs updates of
$\bfSigma$. Sampling from $p(\bZ\mid\bfSigma,\bY)$ is similar to that
in \eqref{eq:con:z} for $p(\bZ\mid\bfOmega,\bY)$ with only the
modification of replacing $\omega_{ij}$ with $\sigma_{ij}.$

\section{Effectiveness of the new methods}
\subsection{Computational efficiency}
\begin{figure}[htbp]
\centering
\includegraphics{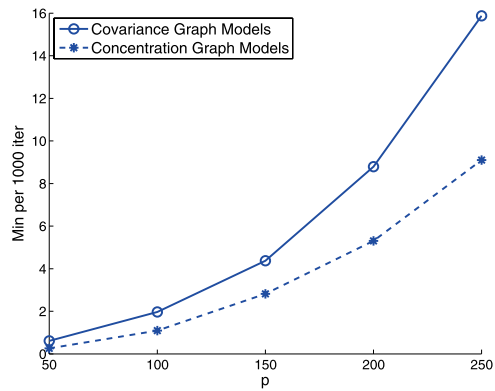}
\caption{Time in minutes for 1000 iterations of SSSL versus dimension
$p$ for covariance graph models (solid) and concentration graph models
(dashed).}
\label{fig:computingtime}
\end{figure}

The computational speed and the scalability of SSSL block Gibbs
samplers are evaluated empirically.
The data of dimension $p\in\{50,100,150,200,250\}$ and sample size
$n=2p$ are first generated from $\N(0,\bI_p)$ and then
standardized. The samplers are implemented under the hyperparameters
$v_0=0.05,h=50,\pi=2/(p-1)$ and $\lambda=1$. All
chains are initialized at the sample covariance matrix. All
computations are implemented on a six-core CPU 3.33GHz desktop using
MATLAB. For each run, we measure the time it takes the block Gibbs
sampler to
sweep across all columns (rows) once, which is called one
iteration. One iteration actually updates each element $a_{ij}$ twice:
once when updating column $i$ and again when updating column $j$. This
property improves its efficiency. The solid and dashed curves in Figure
\ref{fig:computingtime} display the minutes taken for 1000 iterations
versus $p$ for covariance graph models and concentration graph models
respectively.

Overall, the SSSL algorithms run fast. Covariance graph models take
approximately 2 and 9 minutes to generate 1000 iterations for $p=100$
and 200; concentration graph models take even less time, approximately
1.2 and 5 minutes. The relatively slower speed of covariance graph
models is due to a few more matrix inversion steps in updating the
columns in $\bfSigma$. We also measure the mixing of the MCMC output by
calculating the inefficiency factor $1+2\sum_{k=1}^\infty\rho(k)$
where $\rho(k)$ is the sample autocorrelation at lag $k$. We use 5000
samples after 2000 burn-ins
and $K$ lags in the estimation of the inefficiency factors, where $K =
\textrm{argmin}_k\{\rho(k)<2/\sqrt{M},k\geq1 \}$ with $M=5000$ being
the total number of saved
iterations. The median
inefficiency factor among all of the elements of $\bfOmega$ was 1 when $p=100$,
further suggesting the efficiency of SSSL. In our experience, a MCMC
sample of 5000 iterations after 2000 burnins usually generates reliable
results in terms of Monte Carlo errors for $p=100$ node problems,
meaning that the computing time is usually less than 10 minutes, far
less than the few days of computing time required by the existing methods.

%

\subsection{Structure learning accuracy}
The preceding section shows tremendous computational gains of SSSL over
existing methods. This section evaluates these methods on their
structure learning performance using two synthetic scenarios, both of
which are motivated by real-world applications.

\textit{Scenario 1.} The first scenario mimics the dependence pattern
of daily currency exchange rate returns, which has been previously
analyzed via concentration graph models \citep
{CarvalhoWest2006BA,WangReesonCarvalho09}. We use two different
synthetic datasets for the two types of graphical models. For
concentration graph models, \citet{Jones05} generated a simulated
dataset of $p=12$ and $n=250$ that mimics the exchange rate return
pattern. We use their original dataset downloaded from their website.
The true underlying concentration graph has 13 edges and is given by
Figure 2 of \citet{Jones05}. For covariance graph models, we estimate a
sparse covariance matrix with 13 edges based on the data of \citet
{Jones05} and then use this sparse covariance matrix to generate
$n=250$ samples. The true sparse covariance matrix is as follows:
\begin{displaymath}
\left( \scriptsize
\begin{array}{rrrrrrrrrrrr}
$0.239$ & $0.117$ & & & & & & $0.031$ & & & & \\
$0.117$ & $1.554$ & & & & & & & & & & \\
& & $0.362$ & $0.002$ & & & & & & & & \\
& & $0.002$ & $0.199$ & $0.094$ & & & & & & & \\
& & & $0.094$ & $0.349$ & & & & & & & $-0.036$ \\
& & & & & $0.295$ & $-0.229$ & $0.002$ & & & & \\
& & & & & $-0.229$ & $0.715$ & & & & & \\
$0.031$ & & & & & $0.002$ & & $0.164$ & $0.112$ & $-0.028$ & $-0.008$ &
\\
& & & & & & & $0.112$ & $0.518$ & $-0.193$ & $-0.090$ & \\
& & & & & & & $-0.028$ & $-0.193$ & $0.379$ & $0.167$ & \\
& & & & & & & $-0.008$ & $-0.090$ & $0.167$ & $0.159$ & \\
& & & & $-0.036$ & & & & & & & $0.207$
\end{array}
\right).
\end{displaymath}

To assess the performance of structure learning, we compute the number
of true positives (TP) and the number of false positives (FP). We begin
by evaluating some benchmark models against which to compare SSSL. For
concentration graph models, we consider the classical adaptive
graphical lasso \citep{FanFengWu2009} and the Bayesian \iG-Wishart
prior $\W_G(3,\bI_p)$ in (\ref{eq:GW}). For covariance graph models, we
consider the classical adaptive thresholding \citep{CaiLiu2011} and the
Bayesian \iG-inverse Wishart prior $\IW_{G}(3,\bI_p)$ in (\ref
{eq:GIW}). The two classical methods require some tuning parameters,
which are chosen by 10-fold cross-validation. The adaptive graphical
lasso has $\textrm{TP}=8$ and $\textrm{FP}=7$ for concentration graph
models and the adaptive thresholding has $\textrm{TP}=13$ and $\textrm
{FP}=45$ for covariance graph models. They seem to favor less sparse
models, most likely because the sample size is large relative to the
dimensions. The two Bayesian methods are implemented under the priors
of graphs (\ref{eq:ind:bern}) with $\pi=2/(p-1)$. They perform better
than the classical methods. The \iG-Wishart has $\textrm{TP}=9$ and
$\textrm{FP}=1$ for concentration graph models, and the \iG-inverse
Wishart has $\textrm{TP}=7$ and $\textrm{FP}=0$ for covariance graph models.

We investigate the performance of SSSL by considering a range of
hyperparameter values that represent different prior beliefs about
graphs: $\pi\in\{2/(p-1),4/(p-1),0.5\}$, $h \in\{10, 50, 100\}$, and
$v_0 \in\{0.02,0.05,0.1\}$. Under each hyperparameter setting, a
sample of 10000 iterations is used to estimate the posterior median
graph, which is then compared against the true underlying graph.
Table~\ref{tab:con:p12} and \ref{tab:cov:p12} summarize TP and FP for
concentration and covariance graphs respectively. Within each table,
patterns can be observed by comparing results across different values
of one hyperparameter while fixing the others. For $v_0$, a larger
value lowers both TP and FP because $v_0$ is positively related to the
threshold of practical significance that $a_{ij}$ can be treated as
zero. For $\pi$, a larger value increases both TP and FP, and the case
of $\pi=0.5$ seems to produce graphs that are too dense and have high
FP, especially for the concentration graph models in Table~\ref
{tab:con:p12}. For $h$, increasing $h$ reduces both TP and FP, partly
because $h$ is positively related to the threshold that $a_{ij}$ can be
practically treated as zero and in part because $h$ is negatively
related to the implied prior edge-inclusion probability $\Pr(z_{ij})$
as discussed in Panel (a) of Figure \ref{fig:prior:z}. Next, comparing
the results in Tables \ref{tab:con:p12}--\ref{tab:cov:p12} with the two
Bayesian benchmarks reported above, we can see that SSSL is
competitive. For example, Table \ref{tab:cov:p12} shows that, except
for the extreme cases of $\pi=0.5$ that favor graphs that are too dense
or cases of $v_0=0.1$ that favor graphs that are too sparse, SSSL has
approximately the same $\textrm{TP}=7$ and $\textrm{FP}=0$ as the \iG
-inverse Wishart method for covariance graph models.

%
\begin{table}[htbp]
\centering
\scriptsize
\tabcolsep=5.9pt
\caption{Summary of performance measures under different
hyperparameters for a $p=12$ concentration graph example. As for benchmarks,
the classical adaptive graphical lasso has TP=8 and FP=7; the Bayesian
G-Wishart has TP=9 and FP=1. }

\begin{tabular}{lrrrrrrrrrrr}
\toprule
& \multicolumn{3}{c}{$\pi= 2/(p-1)$} & & \multicolumn{3}{c}{$\pi=
4/(p-1)$}& & \multicolumn{3}{c}{$\pi= 0.5$} \\
\cline{2-4} \cline{6-8} \cline{10-12}
$v_0$ & $h=10$ & $h=50$ & $h=100$ & & $h=10$ & $h=50$ & $h=100$& &
$h=10$ & $h=50$ & $h=100$ \\
\midrule
& \multicolumn{11}{c}{\textit{Number of true positives (TP)}} \\
$0.02$ & 9 & 9 & 9 & & 10 & 10 & 10 & & 10 & 10 & 10 \\
$0.05$ & 10 & 9 & 9 & & 10 & 10 & 9 & & 10 & 10 & 10 \\
$0.1$ & 9 & 8 & 8 & & 10 & 8 & 8 & & 10 & 9 & 8 \\
& \multicolumn{11}{c}{\textit{Number of false positives (FP)}} \\
$0.02$ & 3 & 2 & 0 & & 7 & 3 & 1 & & 14 & 4 & 3 \\
$0.05$ & 2 & 0 & 0 & & 4 & 1 & 0 & & 8 & 2 & 0 \\
$0.1$ & 1 & 0 & 0 & & 1 & 0 & 0 & & 4 & 0 & 0 \\
\bottomrule
\end{tabular}

\label{tab:con:p12}
\end{table}
%

%
\begin{table}[htbp]
\centering
\scriptsize
\tabcolsep=5.9pt
\caption{ Summary of performance measures under different
hyperparameters for a $p=12$ covariance graph example. As for benchmarks,
the classical adaptive thresholding has TP=13 and FP=45; the Bayesian
\textit{G}-inverse Wishart has TP=7 and FP=0.}
\begin{tabular}{lrrrrrrrrrrr}
\toprule
& \multicolumn{3}{c}{$\pi= 2/(p-1)$} & & \multicolumn{3}{c}{$\pi=
4/(p-1)$}& & \multicolumn{3}{c}{$\pi= 0.5$} \\
\cline{2-4} \cline{6-8} \cline{10-12}
$v_0$ & $h=10$ & $h=50$ & $h=100$ & & $h=10$ & $h=50$ & $h=100$& &
$h=10$ & $h=50$ & $h=100$ \\
\midrule
& \multicolumn{11}{c}{\textit{Number of true positives (TP)}} \\
$0.02$ & 7 & 7 & 7 & & 8 & 7 & 7 & & 9 & 7 & 7 \\
$0.05$ & 7 & 7 & 6 & & 7 & 7 & 7 & & 7 & 7 & 7 \\
$0.1$ & 5 & 3 & 3 & & 6 & 5 & 4 & & 7 & 5 & 5 \\
& \multicolumn{11}{c}{\textit{Number of false positives (FP)}} \\
$0.02$ & 0 & 0 & 0 & & 1 & 0 & 0 & & 5 & 0 & 0 \\
$0.05$ & 0 & 0 & 0 & & 0 & 0 & 0 & & 0 & 0 & 0 \\
$0.1$ & 0 & 0 & 0 & & 0 & 0 & 0 & & 0 & 0 & 0 \\ \bottomrule
\end{tabular}
\label{tab:cov:p12}
\end{table}

\begin{figure}[htbp]
\centering
\includegraphics{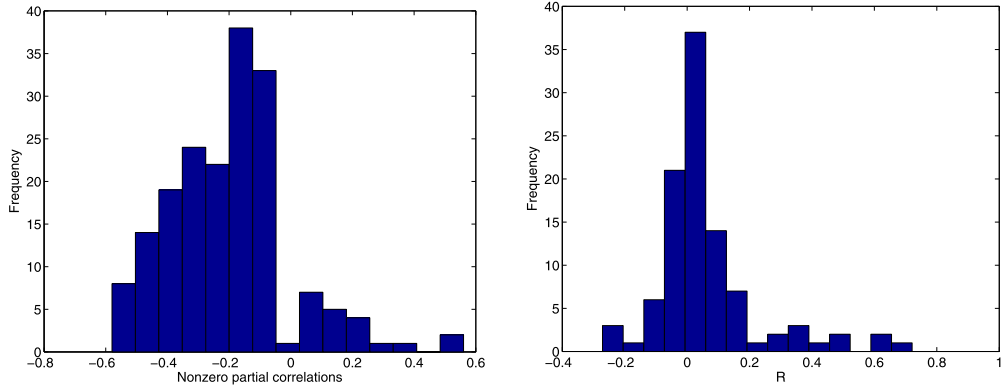}
\caption{Histograms showing the empirical frequency of the non-zero
elements of the partial correlation matrix for the first dataset (left)
and of the nonzero elements of the correlation matrix for the second
dataset (right) in Scenario 2 of $p=150$.}
\label{fig:p150:hist}
\end{figure}

\textit{Scenario 2. } The second scenario mimics the dependence pattern
of gene expression data, for which graphical models are used
extensively to understand the underlying biological relationships. The
real data are the breast cancer data \citep{Jones05,CasteloRoverato06},
which contain $p = 150$ genes related to the estrogen receptor pathway.
Similar to the first scenario, we generate two synthetic datasets for
the two graphical models. For concentration graph models, we first
estimate a sparse concentration matrix with 179 edges based on the real
data, and then generate a sample of 200 observations from this
estimated sparse concentration matrix. For covariance graph models, we
estimate a sparse covariance matrix with 101 edges based on the real
data and then use this sparse covariance matrix to generate a synthetic
data of $n=200$ samples. Panels (a) and (b) of Figure \ref
{fig:p150:hist} display the frequencies of the non-zero partial
correlations and correlations implied by these two underlying sparse
matrices, respectively. Among the nonzero elements, $13\%$ of the
partial correlations and $60\%$ of the correlations are within $0.1$.

%

We repeat the same procedures of fitting the benchmark and proposed
models as in Scenario 1. As for benchmarks, the adaptive graphical
lasso has $\textrm{TP}=145$ and $\textrm{FP}=929$ for concentration
graph models, and the adaptive thresholding has $\textrm{TP}=14$ and
$\textrm{FP}=12$ for covariance graph models. The \iG-Wishart prior has
$\textrm{TP}=105$ and $\textrm{FP}=2$ and takes about four days to run.
The evaluation of the \iG-inverse Wishart is worth elaborating since
existing experiments are conducted in smaller $p$ settings and little
is known about its performance in high-dimensional problems.

The original model fitting algorithm for the \iG-inverse Wishart relies
on the computationally expensive importance sampling for approximating
the normalizing constant and thus is slow and numerically instable for
this $p$=150 dataset. \citet{Silva2013mcmc} develops a new
approximation that requires no Monte Carlo integration, which greatly
speeds up the computation. The MATLAB routine implementing his
algorithm is publicly available on that paper's website. We adopt these
functions with one modification that sets the edge-inclusion
probability to be $2/(p-1)$. The algorithm takes about 2 hours to
complete 1000 sweeps, which is substantially slower than SSSL that
costs about 5 minutes -- see Figure \ref{fig:computingtime}. Since both
SSSL and \citet{Silva2013mcmc} require no Monte Carlo approximation,
the difference in run time is a result of the fact that the SSSL's
block update of $\bZ$ is faster than \citet{Silva2013mcmc}'s
one-edge-at-a-time update.

Using the posterior median graph as an estimate of $G$, the \iG-inverse
Wishart prior $\IW_G(3,\bI)$ produces $\textrm{TP}=3$ and $\textrm
{FP}=3$. These two numbers are surprisingly small. An exploration of
the \iG-inverse Wishart distribution provides some explanations. The
main reason is that when $G$ is sparse and $p$ is large, the \iG
-inverse Wishart prior inadvertently enforces the free elements in
$\bfSigma$ towards zero and hence exerts a strong prior influence on
the posterior distribution. To see this, suppose $G$ is empty, then
(\ref{eq:GIW}) implies that the diagonal elements $\{\sigma_{ii}\}$
follow independent inverse Gamma distributions
\[
p(\sigma_{ii} \mid b) \propto\sigma_{ii}^{-{b+2p\over2}} \exp\{
-{d_{ii} \over2 \sigma_{ii} }\},\quad i=1,\ldots,p,\vadjust{\eject}
\]
which clearly depend on the dimension $p$ and converge to zero rapidly
as $p$ increases for a fixed $b$. Now suppose $G$ is arbitrarily
sparse. Although the theoretical marginal distributions of the free
elements in $\bfSigma$ are unknown, the distribution of $\{\sigma_{ii}\}
$ under an empty graph leads us to conjecture that the free elements in
$\bfSigma$ could be extremely concentrated around zero for large $p$ as
well. In our $p=150$ example, a simulation from $\IW_G(3,\bI_p)$ under
the ground truth $G$ that contains 101 edges supports this conjecture.
The estimated prior mean of these 101 off-diagonal free elements $\{
\sigma_{ij}\}$ is in the range of $ -6 \times10^{-6} $ and $6 \times
10^{-6}$; the estimated prior standard deviation is between $1.9 \times
10^{-4} $ and $2.3 \times10^{-4} $. Such a tightly concentrated prior
provides little probability support for the true graph. The implication
on structure learning is that the Bayes factor might not truly respond
to the data, but largely reflect prior prejudices that $\sigma_{ij}$
are extremely small. In other words, the overly concentrated prior does
not allow the data to speak about $G$ and consequently the posterior
distribution of $G$ is dominated by its prior. In fact, the posterior
sample mean and standard deviation of the number of edges in $\bZ$
computed from the MCMC output of \citet{Silva2013mcmc}
are $170$ and $12.9$, which are close to the prior expected number of
edges, computed as $ {p \choose2} \times{2 \over p-1} =150$ and its
standard deviation, computed as $ \sqrt{ {p \choose2} \times{2 \over
p-1}\times(1-{2 \over p-1}) }=12.16$.

The fundamental cause of this strong prior influence is perhaps that
the parameter space $\{b: b>0\}$ assumed by \cite{SilvaGhahramani2009}
is too restrictive. It might be reasonable to let the parameter space
depend on $G$. For example, the standard inverse-Wishart theory implies
that the parameter space should be $\{b:b>2-2p\}$ when $G$ is empty and
so $b$ could be even negative, and $\{b:b>0\}$ when $G$ is full.
Hierarchical models that allow the value of $b$ to be $G$-dependent
might be helpful. A thorough examination along these lines is beyond
the scope of the current paper. However, it is probably safe to
conclude that further investigation should be called upon to follow the
innovative framework of \cite{SilvaGhahramani2009}.

As for SSSL, Tables \ref{tab:con:p150} and \ref{tab:cov:p150} summarize
its performance. When the results are compared across different levels
of one hyperparameter, the general relations between TP or FP and a
hyperparameter are similar to those in Scenario 1. In fact, the
patterns appear to be more significant in Scenario 2 because priors
have greater influences in this relatively small sample size problem.
When compared with benchmarks, SSSL is competitive. For concentration
graph models, Table \ref{tab:con:p150} suggests that, except for a few
extreme priors that favor overly dense graphs, SSSL produces much
sparser graphs than the classical adaptive graphical lasso, for which
$\textrm{FP}=929$ is too high. When $v_0=0.02$, SSSL is also comparable
to the Bayesian \iG-Wishart prior. When $v_0$ increases, TP drops
quickly because many signals are weak (Figure \ref{fig:p150:hist}) and
are thus treated as practically insignificant by SSSL. For covariance
graph models, Table \ref{tab:cov:p150} suggests that SSSL generally
performs better than the adaptive thresholding. The only exceptions are
cases in which the hyperparameters favor overly dense graphs (e.g., $\pi
=0.5$) or overly sparse graphs (e.g., $v_0=0.1$).
\\

In summary, under sensible choices of hyperparameters, such as
$v_0=0.02,h=50$, and $\pi=2/(p-1)$, SSSL is comparable to existing
Bayesian methods in terms of structure learning accuracy. However,
SSSL's computational advantage of sheer speed and simplicity makes it
very attractive for routine uses.


%
\begin{table}[t!]
\centering
\scriptsize
\tabcolsep=5.9pt
\caption{ Summary of performance measures under different
hyperparameters for a $p=150$ concentration graph example. As for benchmarks,
the classical adaptive lasso has TP=145 and FP=929; the Bayesian
G-Wishart has TP=105 and FP=2. }
\begin{tabular}{lrrrrrrrrrrr}
\toprule
& \multicolumn{3}{c}{$\pi= 2/(p-1)$} & & \multicolumn{3}{c}{$\pi=
4/(p-1)$}& & \multicolumn{3}{c}{$\pi= 0.5$} \\
\cline{2-4} \cline{6-8} \cline{10-12}
$v_0$ & $h=10$ & $h=50$ & $h=100$ & & $h=10$ & $h=50$ & $h=100$& &
$h=10$ & $h=50$ & $h=100$ \\
\midrule
& \multicolumn{11}{c}{\textit{Number of true positives (TP)}} \\
$0.02$ & 106 & 101 & 100 & & 110 & 105 & 102 & & 162 & 136 & 125 \\
$0.05$ & 90 & 82 & 78 & & 97 & 89 & 83 & & 146 & 117 & 109 \\
$0.1$ & 72 & 59 & 56 & & 78 & 67 & 58 & & 122 & 101 & 93 \\
& \multicolumn{11}{c}{\textit{Number of false positives (FP)}} \\
$0.02$& 3 & 0 & 0 & & 9 & 1 & 0 & & 1533 & 238 & 91 \\
$0.05$ & 0 & 0 & 0 & & 0 & 0 & 0 & & 667 & 39 & 9 \\
$0.1$ & 0 & 0 & 0 & & 0 & 0 & 0 & & 169 & 4 & 0 \\
\bottomrule
\end{tabular}
\label{tab:con:p150}%
\end{table}
%

%
\begin{table}[t!]
\centering
\scriptsize
\tabcolsep=5.9pt
\caption{ Summary of performance measures under different
hyperparameters for a $p=150$ covariance graph example. As for benchmarks,
the classical adaptive thresholding has TP=14 and FP=12; the Bayesian
\textit{G}-inverse Wishart cannot run.}
\begin{tabular}{lrrrrrrrrrrr}
\toprule
& \multicolumn{3}{c}{$\pi= 2/(p-1)$} & & \multicolumn{3}{c}{$\pi=
4/(p-1)$}& & \multicolumn{3}{c}{$\pi= 0.5$} \\
\cline{2-4} \cline{6-8} \cline{10-12}
$v_0$ & $h=10$ & $h=50$ & $h=100$ & & $h=10$ & $h=50$ & $h=100$& &
$h=10$ & $h=50$ & $h=100$ \\
\midrule
& \multicolumn{11}{c}{\textit{Number of true positives (TP)}} \\
$0.02$ & 20 & 19 & 15 & & 20 & 20 & 17 & & 38 & 27 & 24 \\
$0.05$ & 12 & 10 & 10 & & 14 & 11 & 10 & & 27 & 15 & 15 \\
$0.1$ & 7 & 6 & 5 & & 7 & 7 & 6 & & 15 & 11 & 10 \\
& \multicolumn{11}{c}{\textit{Number of false positives (FP)}} \\
$0.02$ & 4 & 1 & 1 & & 6 & 3 & 1 & & 1084 & 163 & 69 \\
$0.05$ & 0 & 0 & 0 & & 0 & 0 & 0 & & 229 & 17 & 3 \\
$0.1$ & 0 & 0 & 0 & & 0 & 0 & 0 & & 21 & 0 & 0 \\
\bottomrule
\end{tabular}
\label{tab:cov:p150}%
\end{table}
%
\section{Graphs for credit default swap data}
This section illustrates the practical utility of graphical models by
applying them to credit default swap (CDS) data. CDS is a credit
protection contract in which the buyer of the protection periodically
pays a small amount of money, known as ``spread'', to the seller of
protection in exchange for the seller's payoff to the buyer if the
reference entity defaults on its obligation. The spread depends on the
creditworthiness of the reference entity and thus can be used to
monitor how the market views the credit risk of the reference entity.
The aim of this statistical analysis is to understand the
cross-sectional dependence structure of CDS series and thus the joint
credit risks of the entities. The interconnectedness of credit risks is
important, as it is an example of the systemic risk -- the risk that
many financial institutions fail together.

The data are provided by Markit via Wharton Research Data Services
(WRDS) and comprise daily closing CDS spread quotes for 79
widely traded North American reference entities from January 2, 2001
through April 23, 2012. The quotes are for five-year maturity, as
five-year maturity is the most widely traded and studied term. The
spread quote is then transformed into log returns for graphical model analysis.

To assess the variation in the graphs over time, we estimate graphs
using a one-year moving window. In particular, at the end of each month
$t$, we use the daily CDS returns over the period of month $t-11$ to
month $t$ to estimate the graph $G_t$. The choice of a one-year window
is intended to balance the number of observations, as well as to
accommodate the time-varying nature of the graphs. The first estimation
period begins in January 2001 and continues through December 2001, and
the last is from May 2011 to April 2012. In total, there are 89
estimation periods, corresponding to 89 time-varying graphs for each
type of graph. We set the prior hyperparameters at $v_0=0.02, h=50, \pi
=2/(p-1)$ and $\lambda=1.$ The MCMC are run for 10000 iterations and
the first 3000 iterations are discarded. The graphs are estimated by
the posterior median graph.

%
\begin{figure}[htbp]
\centering
\includegraphics{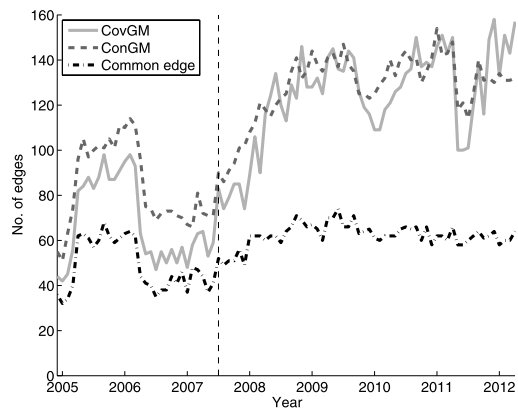}
\caption{Time series plots of the estimated numbers of edges for the
covariance graph (solid line), the concentration graph (dashed), and
the common graph (dash dotted). }
\label{fig:nedge:ts}
\end{figure}

Figure \ref{fig:nedge:ts} shows changes over time of the estimated
number of edges for the two types of graphs and the number of common
edges. The numbers of edges are in the range of 40 and 160 out of a
total of 3081 possible edges, indicating a very high level of sparsity.
At each time point, both types of graph reflect about the same level of
sparsity. Over time, they show similar patterns of temporal variations.
There is a steady upward trend in the number of edges starting from mid
2007 for both types of graph. The timing of the rise of the number of
edges is suggestive. Mid-2007 saw the start of the subprime-mortgage
crisis, which was signified by a number of credit events, including
bankruptcy and the significant loss of several banks and mortgage
firms, such as Bear Stearns, GM finance, UBS, HSBC and Countrywide. If
the number of edges can be viewed as the ``degree of connectedness''
among CDS series, then this observed increase implies that the market
tends to be more integrated during periods of credit crisis and
consequently tends to have a higher systemic risk.

\begin{figure}[t!]
\centering
\includegraphics{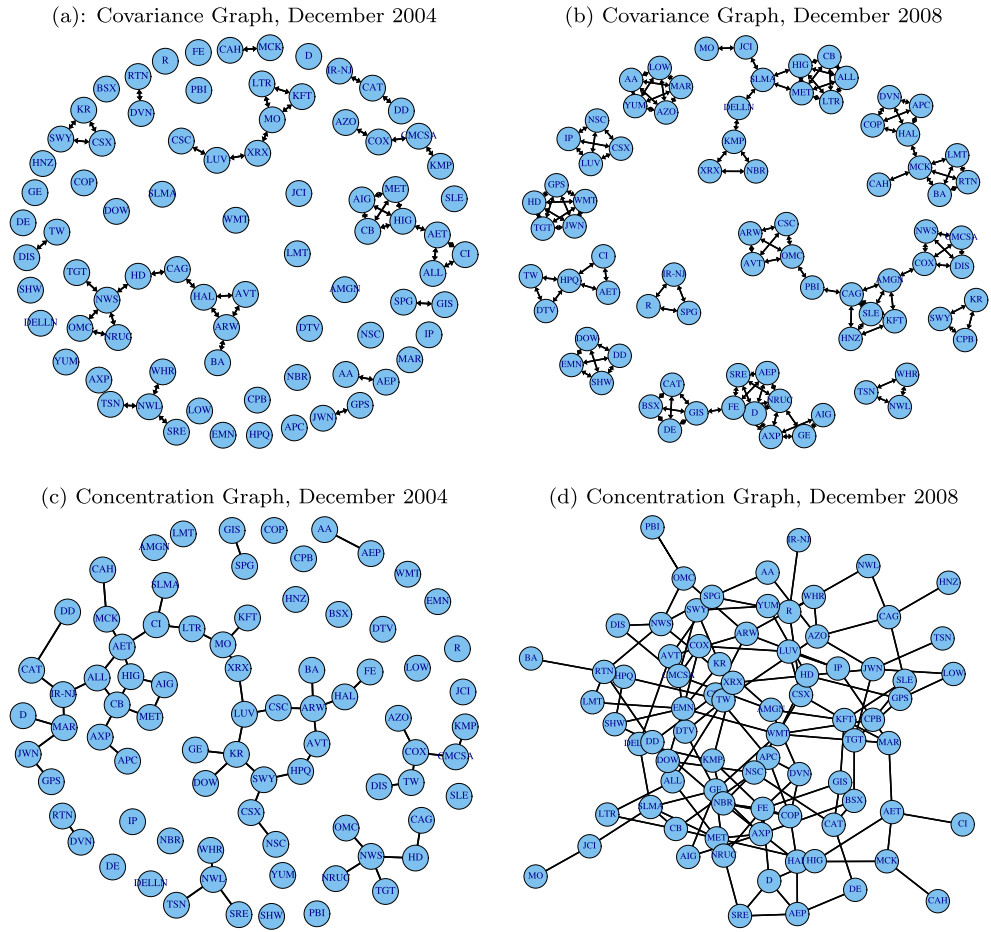}
\caption{Estimated graphs of 79 CDS returns at two different time
points. The four panels are: Panel (a), covariance graph in December
2004; Panel (b), covariance graph in December 2008; Panel (c),
concentration graph in December 2004; and Panel (d) concentration graph
in December 2008.}
\label{fig:graph:snapshop}
\end{figure}

To further illustrate the interpretability of graphs, Figure \ref
{fig:graph:snapshop} provides four snapshots on graph details at two
time points, in December 2004 and in December 2008. The covariance
graph for December 2004 (Panel a) has 44 edges and 30 completely
isolated nodes, as opposed to the covariance graph for December 2008
(Panel b), which has 128 edges and no isolated nodes. These differences
suggest that the connectedness of the network rises as the credit risks
of many reference entities become linked with each other. The same
message is further confirmed by concentration graphs. The concentration
graph for December 2004 (Panel c) has 55 edges and 20 completely
isolated nodes, while the concentration graph for December 2008 (Panel
d) has 135 edges and no isolated nodes. The increase of connectedness
is also manifested by the fact that every pair of nodes in the
concentration graph on December 2008 is connected by a path.

%
%

\begin{figure}[htbp]
\centering
\includegraphics{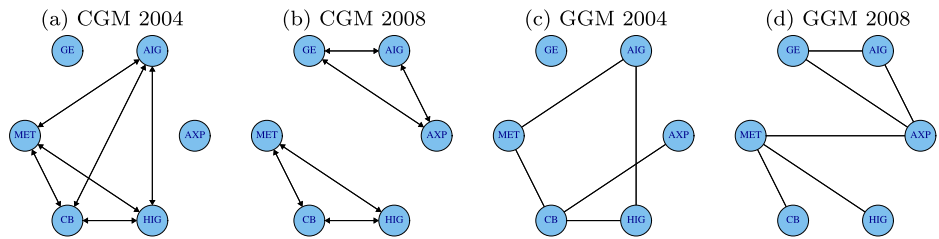}
\caption{ Subgraphs involving AIG estimated for two different time
periods. There are four snapshots: Panel (a), covariance graph in
December 2004; Panel (b), covariance graph in December 2008; Panel (c),
concentration graph in December 2004; and Panel (d) concentration graph
in December 2008. The six subgraph nodes are: American International
Group (AIG); GE Capital (GE); Metropolitan Life Insurance Company
(MET); The Chubb Corporation (CB); Hartford Financial Services Group
(HIG); and American Express Company (AXP).}
\label{fig:graph:snapshop:aig}
\end{figure}

Finally, we zoom into subgraphs involving the American International
Group (AIG) to study whether the graphs make economic sense at the firm
level. AIG provides a good example for study because it suffered a
significant credit deterioration during the 2007--2008 crisis when its
credit ratings were downgraded below ``AA'' levels in September 2008.
Panels (a) and (b) of Figure \ref{fig:graph:snapshop:aig} show the
covariance subgraph in December 2004 and in December 2008. The
subgraphs show that AIG experienced some interesting dependence
structure shifts between these two periods. In December 2004, AIG
formed a clique with three other major insurance companies: the
Metropolitan Life Insurance Company (MET), the Chubb Corporation (CB)
and the Hartford Financial Services Group (HIG). The credit risks of
these four insurance companies are all linked to each other. By
December 2008, the linkages between AIG and the other three insurance
companies disappear, while the linkages among the other three firms
remain. On the other hand, AIG is now connected with two non-insurance
financial companies, the GE Capital (GE) and the American Express
Company (AXP). The concentration subgraphs of AIG displayed in Panels
(c) and (d) of Figure \ref{fig:graph:snapshop:aig} show similar
structural shift patterns as in the covariance graphs, although here
dependence functions differently as conditional dependence. AIG
initially formed a prime component with the other three insurance firms
in December 2004. The connections between AIG and the other insurance
firms were severed in December 2008, and new connections between AIG
and GE and between AIG and AXP arose.

Given that AIG suffered a more severe credit crisis than other
insurance companies, the uncovered network shift might indicate that
the CDS market was able to adjust its view regarding the credit risk of
AIG and treat it as unrelated to its peer insurance firms during the
credit crisis. Such timely changes in dependence structures may shed
light on the question of the information efficiency of the CDS market,
from a point of view that is different from the usual analysis of
changes in levels.

%

\section{Conclusion}
This paper proposes a new framework for the structure learning of
concentration graph models and covariance graph models. The main goal
is to overcome the scalability limitation of existing methods without
sacrificing structure learning accuracy. The key idea is to use
absolutely continuous spike and slab priors instead of the popular
point-mass mixture priors to enable accurate, fast, and simple
structure learning. Our analysis suggests that the accuracy of these
new methods is comparable to that of existing Bayesian methods, but the
model-fitting is vastly faster to run and simpler to implement.
Problems with 100--250 nodes can now be fitted in a few minutes, as
opposed to a few days or even numeric infeasibility under existing
methods. This remarkable efficiency will facilitate the application of
Bayesian graphical models to large problems and will provide scalable
and effective modular structures for more complicated models.

The focus of the paper is on the structure learning of graphs. A
related yet different question is the parameter estimation of the
covariance matrix. Since our priors place zero probability mass on any
sparse matrix containing exact zeros, as opposed to the point-mass
mixture priors, one concern is then about where and how the posterior
of $\bfOmega$ or $\bfSigma$ will concentrate when the true
covariance/concentration matrix is sparse. Our limited experiments
suggest that the posterior distribution from the proposed two-component
normals is indeed more dispersed around zero than those from the
point-mass mixture priors. Take the concentration graph in Scenario 1
as an example. Under SSSL, the average magnitude of the posterior mean
estimates of these $\{\omega_{ij}\}$ corresponding to the true zeros in
$\bfOmega$ is in the range of 0.0085 and 0.045, depending on the
hyperparameters of $(v_0,h,\pi)$. In contrast, under the \iG-Wishart
prior, the estimates of these $\{\omega_{ij}\}$ have a smaller average
magnitude of 0.0040. The small-variance normal shrinks parameters less
aggressively than the point-mass mixture. Although we do not find it
problematic for structure learning, such weaker shrinkage may cause the
SSSL's performance of parameter estimation to be suboptimal. A
refinement is to replace the normal distribution with distributions
having higher densities near zero. The Bayesian shrinkage regression
literature shows that some heavy-tailed distributions offer comparable
parameter estimation performance to the point-mass mixture priors
(e.g., \citealt{Armagan11,GriffinBrown10,CarvalhoPolsonScott2010}).
Computational tractability of SSSL is maintained by applying the
data-augmentation to the mixture normal representation of these
alternative distributions.

\section*{Appendix}
\subsection*{Proof of Proposition \ref{propsn1}}
Clearly, the conditional distribution of $\bfOmega$ given the
edge-inclusion indicators $\bZ$ is
\begin{multline}\label{eq:jointpost:con}
p(\bfOmega \mid\bY,\bZ) \propto|\bfOmega|^{{n\over2}}\exp\{-\tr
({1\over2} \bS\bfOmega) \} \prod_{i<j} \bigg\{ \exp(-{\omega
^2_{ij}\over2v^2_{z_{ij}}}) \bigg\} \prod_{i=1}^p \bigg\{\exp
(-{{\lambda\over2}} \omega_{ii}) \bigg\}.
\end{multline}
Under the partitions \eqref{eq:partition:con}, the conditional
distribution of the last column in $\bfOmega$ is
\begin{multline}
p(\bfomega_{12},\omega_{22} \mid\bY,\bZ,\bfOmega_{11})\propto(\omega
_{22} - \bfomega_{12}^\prime\bfOmega_{11}^{-1} \bfomega_{12})^{{n\over
2}} \\ \exp\big[- {1\over2} \{ \bfomega_{12}^\prime\bD^{-1} \bfomega
_{12} + 2 \bs_{12}^\prime\bfomega_{12}+(s_{22}+\lambda) \omega_{22}\}
\big], \nonumber
\end{multline}
where $\bD= \diag(\bv_{12})$. Consider a change of variables $(\bu
,\omega_{22}) \rightarrow(\bu=\bfomega_{12}, v = \omega_{22} - \bfomega
_{12}^\prime\bfOmega_{11}^{-1} \bfomega_{12}),$ whose
Jacobian is a constant not involving $(\bfomega_{12},v)$. So
\[
p(\bu,v \mid\bY,\bZ,\bfOmega_{11})\propto v^{{n\over2}} \exp(
-{s_{22}+\lambda\over2}v)\exp\bigg(- {1\over2} \big[ \bu^\prime\{\bD
^{-1}+(s_{22}+\lambda) \bfOmega_{11}^{-1} \} \bu+ 2 \bs_{12}^\prime\bu
\big] \bigg).
\]
This implies that:
\[
(\bu,v) \mid (\bfOmega_{11},\bZ,\bY) \sim \N(-\bC\bs_{12},\bC)
\textsc{Ga}({n\over2}+1,{s_{22} + \lambda\over2}),
\]
where $\bC= \{(s_{22}+\lambda)\bfOmega_{11}^{-1}+\bD^{-1}\}^{-1}$.

\subsection*{Proof of Proposition \ref{propsn2}}
Given edge-inclusion indicators $\bZ$, the conditional posterior
distribution of $\bfSigma$ is
\begin{align}\label{eq:jointpost:cov}
p(\bfSigma\mid\bY,\bZ) \propto|\bfSigma|^{-{n\over2}}\exp\{-{1\over
2} \tr(\bS\bfSigma^{-1}) \}\prod_{i<j} \bigg\{ \exp(-{\sigma
^2_{ij}\over2v^2_{z_{ij}}}) \bigg\} \prod_{i=1}^p \bigg\{\exp
(-{{\lambda\over2}} \sigma_{ii}) \bigg\}.
\end{align}
%
%
Under partitions (\ref{eq:partition:cov}), consider a transformation
$(\bfsigma_{12},\sigma_{22})\rightarrow(\bu=\bfsigma_{12}, v = \sigma
_{22} - \bfsigma_{12}^\prime\bfSigma_{11}^{-1} \bfsigma_{12}),$ whose
Jacobian is a constant not involving $(\bu,v)$ and apply the block
matrix inversion to $\bfSigma$ using blocks $(\bfSigma_{11},\bu,v)$:
%
\begin{equation}\label{eq:blockinv}
\bfSigma^{-1} = \left(
\begin{array}{cc}
\bfSigma_{11}^{-1}+\bfSigma_{11}^{-1} \bu\bu^\prime\bfSigma
_{11}^{-1}v^{-1} & -\bfSigma_{11}^{-1} \bu v^{-1} \\ -\bu^\prime
\bfSigma_{11}^{-1}v^{-1} & v^{-1}
\end{array}
\right).
\end{equation}
After removing some constants not involving $(\bu,v)$, the terms in
(\ref{eq:jointpost:cov}) can be expressed as a function of $(\bu,v)$:
\begin{eqnarray*}
|\bfSigma| &\propto& v, \nonumber\\
\tr(\bS\bfSigma^{-1}) &\propto& \bu^\prime\bfSigma_{11}^{-1} \bS_{11}
\bfSigma_{11}^{-1} \bu v^{-1}-2\bs_{12}^\prime\bfSigma_{11}^{-1}\bu
v^{-1}+s_{22}v^{-1},
\\ \prod_{i<j} \bigg\{ \exp(-{\sigma^2_{ij}\over2v^2_{z_{ij}}}) \bigg
\} &\propto& \exp(-{1\over2} \bu^\prime\bD^{-1} \bu) , \\
\prod_i \bigg\{\exp(-{{\lambda\over2}} \sigma_{ii}) \bigg\} &\propto&
\exp(-{1\over2}\lambda(\bu^\prime\bfSigma_{11}^{-1} \bu+v) ) ,
\end{eqnarray*}
where $\bD= \diag(\bv_{12})$. Holding all but $(\bu,v)$ fixed, we can
then rewrite the logarithm of (\ref{eq:jointpost:cov}) as\vadjust{\eject}
\begin{align*}
\log\,p(\bu,v \mid- ) = -{1\over2} \biggr\{ n\log(v) + \bu^\prime
\bfSigma_{11}^{-1} \bS_{11} \bfSigma_{11}^{-1} \bu v^{-1}-2\bs
_{12}^\prime\bfSigma_{11}^{-1}\bu v^{-1} + s_{22}v^{-1} \nonumber\\ +
\bu^\prime\bD^{-1} \bu+ \lambda\bu^\prime\bfSigma_{11}^{-1} \bu+
\lambda v \biggr\}+constant.
\end{align*}
This gives the conditionals of $\bu$ and $v$ as
\begin{align*}(\bu\mid\bY,\bZ,\bfSigma_{11},v) & \sim \N\big\{(\bB+
\bD^{-1})^{-1} \bw, (\bB+\bD^{-1})^{-1} \big\}, \\ (v \mid\bY,\bZ
,\bfSigma_{11},\bu) & \sim\GIG(1-n/2,\lambda, \bu^\prime\bfSigma
_{11}^{-1} \bS_{11} \bfSigma_{11}^{-1} \bu-2\bs_{12}^\prime\bfSigma
_{11}^{-1}\bu+ s_{22}),
\end{align*}
where $\bB= \bfSigma_{11}^{-1} \bS_{11} \bfSigma_{11}^{-1}v^{-1}+\lambda
\bfSigma_{11}^{-1}$, $\bw= \bfSigma_{11}^{-1}\bs_{12}v^{-1}$ and $\GIG
(q,a,b)$ denotes the generalized inverse Gaussian distribution with
probability density function:
\[
p(x) = \frac{(a/b)^{q/2}}{2 K_q(\sqrt{ab})} x^{(q-1)} e^{-(ax + b/x)/2},
\]
with $K_q$ a modified Bessel function of the second kind.

\section*{Supplementary materials}
The MATLAB routines implementing all frequentist and Bayesian
procedures used in the paper are available from the author's web
site of the paper at: \texttt{\surl{https://www.msu.\\edu/~haowang/RESEARCH/SSSL/sssl.html}}.

\bibliographystyle{ba}

\begin{acknowledgement}
The author would like to thank two anonymous referees for helpful
comments that improved
the quality of the manuscript.
\end{acknowledgement}

\end{document}